\documentclass{amsart}

\usepackage{amsmath,amssymb,amsfonts}
\usepackage{amsthm}

 \usepackage{listings}
 \usepackage{enumitem}
 \usepackage{accents} 
\usepackage{tikz}
 \usepackage{graphicx}
\usepackage{wrapfig}
\usetikzlibrary{calc,hobby}

\newtheorem{theorem}{Theorem}[section]
\newtheorem{lemma}[theorem]{Lemma}

\newtheorem{proposition}[theorem]{Proposition}

\theoremstyle{definition}
\newtheorem{definition}[theorem]{Definition}
\newtheorem{example}[theorem]{Example}

\theoremstyle{remark}
\newtheorem{remark}[theorem]{Remark}

\numberwithin{equation}{section}

\providecommand{\abs}[1]{\lvert#1\rvert}
\providecommand{\norm}[1]{\lVert#1\rVert}
\begin{document}

\title[Continuity of extensions of Lipschitz maps and of monotone maps]{Continuity of extensions of Lipschitz maps and of monotone maps}


\author{Krzysztof J. Ciosmak}

\address{Department of Mathematics, University of Toronto, Bahen Centre, 40 St. George St., Room 6290, Toronto, Ontario, M5S 2E4, Canada}
 \address{Fields Institute for Research in Mathematical Sciences, 222 College Street, Toronto, Ontario, M5T 3J1, Canada }
\email{k.ciosmak@utoronto.ca}
 \thanks{Part of the research presented in this paper was performed when the author was a post-doctoral fellow at the University of Oxford, supported by the European Research Council Starting Grant CURVATURE, grant agreement No. 802689. 
}



\date{}


\begin{abstract}

Let $X$ be a subset of a Hilbert space. 
We prove that if $v\colon X\to \mathbb{R}^m$ is such that
\begin{equation*}
    \Big\lVert v(x)-\sum_{i=1}^m t_iv(x_i)\Big\rVert\leq     \Big\lVert x-\sum_{i=1}^m t_ix_i\Big\rVert
\end{equation*}
for all $x,x_1,\dotsc,x_m\in\mathbb{R}^m$ and all non-negative $t_1,\dotsc,t_m$ that add up to one, then for any $1$-Lipschitz $u\colon A\to\mathbb{R}^m$, with $A\subset X$, there exists a $1$-Lipschitz extension $\tilde{u}\colon X\to\mathbb{R}^m$ of $u$ such that
the uniform distance on $X$ between $\tilde{u}$ and $v$ is the same as the uniform distance on $A$ between $u$ and $v$.

Moreover, if either $m\in \{1,2,3\}$ or $X$ is convex, we prove the converse: we show that a map $v\colon X\to\mathbb{R}^m$ that allows for a $1$-Lipschitz, uniform distance preserving extension  of any $1$-Lipschitz map on a subset of $X$ also satisfies the above set of inequalities.

We also prove a similar continuity result concerning extensions of monotone maps.

Our results hold true also for maps taking values in infinite-dimensional spaces.
\end{abstract}

\maketitle

\section{Introduction}

Let $Y,Z$ be two Hilbert spaces. Let $X$ be any subset of a $Z$. We say that a map $u\colon X\to Y$ is $1$-Lipschitz if for any $x,y\in X$ we have
\begin{equation*}
\norm{u(x)-u(y)}\leq \norm{x-y}.
\end{equation*}
A theorem of Kirszbraun, see \cite{Kirszbraun1934}, tells that any $1$-Lipschitz map on a subset of a finite-dimensional Hilbert space may be extended to a $1$-Lipschitz map on the entire space. This result has been later generalised  by Valentine, see \cite{Valentine1945}, so that it is applicable to $1$-Lipschitz maps between arbitrary Hilbert spaces

\begin{theorem}\label{thm:Kirszbraun}
Let $Y,Z$ be two Hilbert spaces. Let $X\subset Z$. Let $u\colon X\to Y$ be a $1$-Lipschitz map. Then there exists a $1$-Lipschitz map $\tilde{u}\colon Z\to Y$ such that $\tilde{u}|_{X}=u$.
\end{theorem}

Suppose that $Y,Z$ are two Hilbert spaces and that $X\subset Z$. Let $v\colon X\to Y$ be a map. We shall be concerned with finding conditions on $v$ that are sufficient and necessary for the existence, given a $1$-Lipschitz map $u\colon A\to Y$ on a subset $A\subset X$, of a $1$-Lipschitz extension $\tilde{u}\colon X\to Y$ of $u$ such that
\begin{equation*}
    \sup\big\{\norm{\tilde{u}(x)-v(x)}\mid x\in X\big\}=    \sup\big\{\norm{u(x)-v(x)}\mid x\in A\big\}.
\end{equation*}
Phrased in words, we are interested in whether it is possible to find a $1$-Lipschitz extension of $u$ to $X$ that would be not farther to $v$ on $X$ than $u$ is to $v$ on $A$. The main result that we shall prove is the following theorem.

\begin{theorem}\label{thm:main}
Let $Z,Y$ be two real Hilbert spaces. Suppose that $X\subset Z$ and let $v\colon X\to Y$ be a map. Then each of the following conditions implies the subsequent one:
\begin{enumerate}[label=\roman*)]
    \item\label{i:v} for all finite $m\leq \mathrm{dim}Y$, all $x,x_1,\dotsc,x_m\in X$ and all non-negative real numbers $t_1,\dotsc,t_m$ that sum up to one we have  
\begin{equation*}
\Big\lVert v(x)-\sum_{i=1}^mt_iv(x_i)\Big\rVert\leq \Big\lVert x-\sum_{i=1}^mt_ix_i\Big\rVert,
\end{equation*}
\item\label{i:K} for any closed, convex set $K\subset Y$, for any $A\subset X$, and for any $1$-Lipschitz map $u\colon A\to Y$ such that for all $x\in A$
\begin{equation*}
u(x)-v(x)\in  K
\end{equation*}
there exists a $1$-Lipschitz extension $\tilde{u}\colon X\to Y$ of $u$ such that for all $x\in X$
\begin{equation*}
\tilde{u}(x)-v(x)\in K,
\end{equation*}
\item\label{i:ball} for any $A\subset X$, any $\delta\geq 0$, and for any $1$-Lipschitz map $u\colon A\to Y$ such that for all $x\in A$
\begin{equation*}
\norm{u(x)-v(x)}\leq \delta
\end{equation*}
there exists a $1$-Lipschitz extension $\tilde{u}\colon X\to Y$ of $u$ such that for all $x\in X$
\begin{equation*}
\norm{\tilde{u}(x)-v(x)}\leq\delta.
\end{equation*}
\end{enumerate}
Moreover, if either 
\begin{enumerate}
    \item $\mathrm{dim}Y\in\{1,2,3\}$, or,
    \item $X$ is convex,
\end{enumerate}
then \ref{i:ball} implies \ref{i:v}, so that \ref{i:v}-\ref{i:ball} are equivalent.
\end{theorem}

Let us note that the theorem shows that if $v$ satisfies \ref{i:v}, then we may find a $1$-Lipschitz extension $\tilde{u}\colon X\to Y$ of $u\colon A\to Y$ such that 
\begin{equation*}
    \sup\big\{\norm{\tilde{u}(x)-v(x)}_K\mid x\in X\big\}=    \sup\big\{\norm{u(x)-v(x)}_K\mid x\in A\big\},
\end{equation*}
where $\norm{\cdot}_K$ is a norm on $Y$ of our choice. 
The fact that \ref{i:v} implies \ref{i:ball} has already been proven in \cite{Ciosmak2020}, where it has been also conjectured that \ref{i:v} and \ref{i:ball} are equivalent. Our theorem therefore provides an affirmative answer to this conjecture, under an additional assumption that either $\mathrm{dim}Y\in \{1,2,3\}$ or $X$ is convex.

If $X$ is convex, the theorem above could be readily inferred from \cite[Theorem 3.9]{Ciosmak2021}, which tells that Theorem \ref{thm:main} holds true if $X=Z$ is the entire Hilbert space. 

Let us observe a relation of Theorem \ref{thm:main} to the fundamental theorem of affine geometry, see e.g., \cite[Theorem 5.4.]{Artstein2017}. The theorem says that an injective map $v\colon\mathbb{R}^n\to\mathbb{R}^n$, with $n\geq 2$, that maps any three collinear points onto collinear points, is necessarily affine. Clearly, if $X$ is convex, and if $\mathrm{dim}Y\geq 2$, then any map satisfying \ref{i:v} is affine. Thus, Theorem \ref{thm:main} gives a way to characterise affine and $1$-Lipschitz maps with values in an at least two-dimensional space, that resembles the fundamental theorem of affine geometry. Let us note that both characterisations of affine maps do not hold in the one-dimensional setting, see \cite[Proposition 3.8]{Ciosmak2021}.



\subsection{Motivation -- multi-dimensional localisation}

The motivation for this work stems from considerations pertaining to variational problems concerning $1$-Lipschitz maps, including in particular the optimal transport problem, see \cite{Villani2003, Villani2009} for an extensive account on this topic. As described in detail in \cite{Ciosmak2021}, when the cost function in an optimal transport problem is a metric, then the existence of continuous $1$-Lipschitz extensions of real-valued functions can be used to prove the so-called mass balance condition, which plays a crucial role in the localisation technique. The technique was developed in the setting of Riemannian manifolds by Klartag in \cite{Klartag2017}, and in the setting of metric measure spaces by Cavalletti and Mondino in \cite{Cavalletti20172, Cavalletti2017}, and allowed for proofs of various novel and sharp geometric and functional inequalities. Recent works \cite{Ciosmak2023,Ciosmak20232} showed a relation of the localisation technique to martingale transport.

It was conjectured in \cite[Chapter 6]{Klartag2017} that if we consider in place of real-valued $1$-Lipschitz functions their vector-valued analogues and if we look at the variational problem
\begin{equation*}
\sup\Big\{\int_{\mathbb{R}^n}\langle u,d\mu\rangle\mid u\colon\mathbb{R}^n\to\mathbb{R}^m\text{ is }1\text{-Lipschitz}\Big\},
\end{equation*}
where $\mu$ is a $\mathbb{R}^m$-valued Borel measure such that $\mu(\mathbb{R}^n)=0$, then a multi-dimensional mass balance condition will also be satisfied. That is, let $v\colon\mathbb{R}^n\to\mathbb{R}^m$ be a $1$-Lipschitz map that attains the above supremum. A \emph{leaf} $\mathcal{S}$ of $v$ is a maximal subset of $\mathbb{R}^n$ such that $v|_{\mathcal{S}}$ is an isometry. The conjecture says if $\mu$ is absolutely continuous, then  $\mu(A)=0$ for any Borel set $A$ that is a union of leaves of $v$.
As shown in \cite[Theorem 3.9]{Ciosmak2021}, the argument outlined in \cite{Klartag2017} contains a gap. In \cite[Theorem 5]{Ciosmak20214} it is shown that in fact the conjecture is false if $m>1$. Theorem \ref{thm:main} shows that if $\mu$ is concentrated on a set $X$, and if the maximiser satisfies \ref{i:v} of Theorem \ref{thm:main}, then the conjecture of Klartag holds true. We shall not describe in detail the proof. Instead, we point out that the proof presented in \cite[Corollary 4.5]{Klartag2017} works, mutatis mutandis, in this case. Let us point out work \cite{Ciosmak20212}, where the partitioning of $\mathbb{R}^n$ into leaves of a given $1$-Lipschitz map is studied thoroughly.

We shall note also a relation of the extension problem in the discrete setting to a problem posed by Alon in \cite{Alon1990}. Suppose that we are given $n$ sets $A_1,\dotsc, A_n\subset\mathbb{R}^n$, each of cardinality $d$, whose union is in general position. 
The problem is to find an effective way of computing pairwise disjoint simplices $S_1,\dotsc,S_d$ such that $\abs{S_i\cap A_j}=1$ for each $i=1,\dotsc,d$ and $j=1,\dotsc,n$. Existence of such simplices is ensured by a result of Akiyama and Alon proven in \cite{Akiyama1989}.
If $n=2$, the problem is solved through the optimal transport problem. Namely, it suffices to find a bijection $f\colon A_1\to A_2$ for which the sum of the lengths of all segments 
joining the points $a$ and $f(a)$, $a\in A_1$, is minimal. The triangle inequality implies that the resulting line segments are pairwise disjoint. Another version of the solution of this problem for $n=2$ is to look at the leaves of a Kantorovich potential associated with the transport problem. This shows a link to Theorem \ref{thm:main}. 
If $n>2$, then no solution to the Alon problem is known. 


\subsection{Kirszbraun theorem and related literature}

Among a number of the proofs of the Kirszbraun theorem some, e.g., those presented in \cite{Kirszbraun1934, Valentine1945, Schoenberg1953}, use the Kuratowski--Zorn lemma and thus are non-constructive. 
Other proofs, among them \cite{Akopyan2008, Brehm1981}, are constructive. There exists also a formula for the extension, see \cite{Azagra2021}.
Let us also note the approaches that use the Fenchel duality and Fitzpatrick functions, see \cite{Bauschke2007, Reich2005}. Related methods also allowed for a constructive proof, see \cite{Bauschke2010}. These methods allow for further and novel contributions to the field. Let us mention a recent work \cite{Cavagnari2023}, where an invariant version of the Kirszbraun theorem is proven.

We refer the reader also to the work of Dacorogna and Gangbo \cite{Dacorogna2006} where various extensions properties of vector-valued maps are studied.
In \cite{Pak2019} a version of the Kirszbraun theorem for graphs is provided.

Note that Kirszbraun theorem holds not only in Euclidean spaces, but also for spaces with an upper or lower bound on the curvature in the sense of Alexandrov as proven by Lang and Schroeder in \cite{Lang1997}.

We refer to \cite{Cobzac2019} to a vast survey of results on Lipschitz maps. 
For references concerning extensions problems in other settings we refer to \cite{Ciosmak2021}, which include the problem of Whitney of extending differentiable functions, studies of which were commenced in \cite{Whitney1934}.


\subsection{Continuity of extensions}

In \cite{Ciosmak2021} we studied the rate of continuity of extensions of Lipschitz maps. Namely, we supposed that we were given two sets $A\subset B\subset\mathbb{R}^n$ and $1$-Lipschitz maps $u\colon A\to\mathbb{R}^m$ and $v\colon B\to\mathbb{R}^m$, with $m>1$.
We were interested in 
\begin{equation*}
\inf \big\{\sup\big\{ \norm{\tilde{u}(x)-v(x)}\mid x\in B\big\}\mid\tilde{u}\colon B\to\mathbb{R}^m\text{ is }1\text{-Lipschitz extension of }u\big\}.
\end{equation*}
We have shown that for any $u,v$ this quantity is bounded from above by 
\begin{equation*}
 \sqrt{\delta^2+2\delta d_v(A,B)},
\end{equation*}
where 
\begin{equation*}
    d_v(A,B)=\sup\big\{\norm{v(x)-v(y)}\mid x\in A, y\in B\big\}\text{ and  }\delta=\sup\big\{\norm{v(x)-u(x)}\mid x\in A\big\}.
\end{equation*}
Moreover, as shown in \cite{Ciosmak2021}, this bound is sharp.

In \cite{Kopecka20122, Kopecka2012, Kopecka2011} a related question concerning the existence of a Lipschitz constant preserving continuous extension operator on the space of Lipschitz maps is studied.
\subsection{Monotone maps}

Let $Y$ be a linear topological space, let $Y^*$ denote its dual space, let $X\subset Y$. Let us recall that a map $u\colon X\to Y^*$ is called monotone whenever
\begin{equation*}
    \langle u(x)-u(y),x-y\rangle\geq 0\text{ for all }x,y\in X.
\end{equation*}
Existence of maximal monotone operators is well-established to be related to extensions of Lipschitz maps. The topic of maximal monotone operators has been initiated by Minty, \cite{Minty1961,Minty1962}, studied in \cite{Fitzpatrick1989}, and has since been an active area of research. Let us mention recent works \cite{Cavagnari2023, Ghoussoub2008, Visintin2017} pertaining to the topic. However, to the author's knowledge, no results concerning continuity of monotone extensions are known.

Let $v\colon X\to Y^*$ and let $A\subset X$. We shall study conditions on $v$ under which for any monotone $u\colon A\to Y^*$, there exists a monotone extension $\tilde{u}\colon X\to Y^*$ of $u$ that is no farther to $v$ on $X$ than $u$ is to $v$ on $A$.

We shall prove the following theorem.

\begin{theorem}\label{thm:monotone}
Let $Y$ be a linear topological space. Suppose that $X\subset Y$. Let $v\colon X\to Y^*$. Suppose that for all finite $m\leq \mathrm{dim}(Y^*)$, all $x,x_1,\dotsc,x_m\in X$ and all non-negative real numbers $t_1,\dotsc,t_m$ that sum up to one we have
\begin{equation*}
\Big\langle v(x)-\sum_{i=1}^mt_iv(x_i),x-\sum_{i=1}^mt_ix_i\Big\rangle\geq 0.
\end{equation*}
Then for any closed, bounded, convex set $K\subset Y^*$, for any $A\subset X$, and for any monotone map $u\colon A\to Y^*$ such that for all $x\in A$
\begin{equation*}
u(x)-v(x)\in  K
\end{equation*}
there exists a monotone extension $\tilde{u}\colon X\to Y^*$ of $u$ such that for all $x\in X$
\begin{equation*}
\tilde{u}(x)-v(x)\in K.
\end{equation*}
\end{theorem}

\subsection{Maps of semi-bounded strain}

In the study of Michell trusses, see, e.g., \cite{Michell1904,Michell1899, Michell1900} for foundational works and and \cite{Rindler2023,Gangbo2008, Olbermann2017,Olbereman2022} for recent contributions, a class of functions of bounded strain plays an important role. The class consists of functions $u\colon A\to\mathbb{R}^n$, where $A\subset\mathbb{R}^n$, such that
\begin{equation}\label{eqn:daco}
    \abs{\langle u(x)-u(y),x-y\rangle}\leq \norm{x-y}^2\text{ for all }x,y\in A.
\end{equation}
The extensions properties of this class have been studied by Dacorogna and Gangbo in \cite{Dacorogna2006}. In particular, it is shown in \cite[Theorem 22]{Dacorogna2006} that there exist sets $A\subset\mathbb{R}^n$ and functions $u\colon A\to\mathbb{R}^n$ satisfying (\ref{eqn:daco}) with no extension to $\mathbb{R}^n$.

The situation is different, when one considers $u\colon A\to\mathbb{R}^n$ that satisfy
\begin{equation}\label{eqn:dacoba}
    \langle u(x)-u(y),x-y\rangle\leq \norm{x-y}^2\text{ for all }x,y\in A.
\end{equation}
In \cite[Proposition 28, (i)]{Dacorogna2006}, it is shown that if $A$ is a finite set, then $u$ can be extended to a larger finite set. However, as shown in \cite[Proposition 28, (ii)]{Dacorogna2006}, there are infinite sets $A$ and functions $u$ satisfying (\ref{eqn:dacoba}), with no extension to a larger set.

\begin{definition}
 Let $Y$ be a Hilbert space.  Let $A\subset Y$, $u\colon A\to Y$. Then $u$ is said to be of $1$-semi-bounded strain if 
 \begin{equation*}
    \langle u(x)-u(y),x-y\rangle\leq \norm{x-y}^2\text{ for all }x,y\in A.
\end{equation*}
\end{definition}
Suppose now that $v\colon X\to\mathbb{R}^n$, $A\subset X$. Here we show that, under certain conditions on $v$, not only does there exist an extension of any $u\colon A\to\mathbb{R}^n$ of $1$-semi-bounded strain  to a map $\tilde{u}\colon X\to\mathbb{R}^n$ of $1$-semi-bounded strain on $X$, but also such $\tilde{u}$ can be chosen to be no farther to $v$ on $X$ than $u$ is to $v$ on $A$.

This shows that if we impose certain boundedness properties on $u$, then an extension always exists.

\begin{theorem}\label{thm:michell}
Let $Y$ be a real Hilbert space. Suppose that $X\subset Y$. Let $v\colon X\to Y$. Suppose that for all finite $m\leq \mathrm{dim}Y$, all $x,x_1,\dotsc,x_m\in X$ and all non-negative real numbers $t_1,\dotsc,t_m$ that sum up to one we have
\begin{equation*}
\Big\langle v(x)-\sum_{i=1}^mt_iv(x_i),x-\sum_{i=1}^mt_ix_i\Big\rangle\leq \Big\lVert x-\sum_{i=1}^mt_ix_i\Big\rVert^2.
\end{equation*}
Then for any closed, bounded, convex set $K\subset Y$, for any $A\subset X$, and for any map $u\colon A\to Y$ of $1$-semi-bounded strain such that for all $x\in A$
\begin{equation*}
u(x)-v(x)\in  K
\end{equation*}
there exists an extension $\tilde{u}\colon X\to Y$ of $u$ of $1$-semi-bounded strain such that for all $x\in X$
\begin{equation*}
\tilde{u}(x)-v(x)\in K.
\end{equation*}
\end{theorem}

\begin{remark}
    If we put $v=0$, then Theorem \ref{thm:michell} says that any map of $1$-semi-bounded strain that has bounded range, admits an extension to a map of $1$-semi-bounded strain on the entire Hilbert space.
\end{remark}

The theorem is, to the author's knowledge, the first result that gives a condition under which a map of $1$-semi-bounded strain admits an extension and the first result concerned with continuity of extensions of maps of $1$-semi-bounded strain. 

Note that a map $v\colon X\to Y$ is of $1$-semi-bounded strain if and only if $\mathrm{id}-v\colon X\to Y$ is monotone. Therefore the above result is an immediate consequence of Theorem \ref{thm:monotone}.


\subsection*{Structure of the paper}
In Section \ref{s:minty} we provide a statement of the Minty theorem, the main tool for proving our extension theorems.

In Section \ref{s:existence} we prove part of Theorem \ref{thm:main}. In particular, the section contains a proof that \ref{i:v} implies \ref{i:K} in Theorem \ref{thm:main}.

Section \ref{s:necess} completes  the proof of Theorem \ref{thm:main}, i.e., proves that \ref{i:ball} implies \ref{i:v} for $\mathrm{dim}Y\in \{1,2,3\}$ or when $X$ is convex. Section \ref{ss:leaves} provides auxiliary results, Section \ref{ss:conseq} outlines the main punch line, Section \ref{ss:construction} yields a construction necessary for the punch line.
Section \ref{ss:convex} is concerned with the case when $X$ is convex.

In Section \ref{s:monotone} we provide a proof of Theorem \ref{thm:monotone}. 

In Section \ref{s:further} we discuss further possible generalisations and outline one of them.

\section{Minty theorem}\label{s:minty}

We shall start by recalling one of the main tools that we shall use in the proof of Theorem \ref{thm:main}. It is the theorem of Minty, see \cite{Minty1970}, which encompasses several Kirszbraun-type theorems. We shall first bring to the reader's attention the following definition.

\begin{definition}\label{defin:Ka}
Let $Y$ be a real vector space and $X$ be a set. A real function on $Y$ is called finitely lower semi-continuous if its restriction to any finitely-dimensional subspace of $Y$ is lower semi-continuous. A function $\Phi\colon Y\times X\times X\to\mathbb{R}$ is called a Kirszbraun function if $\Phi$ is both finitely lower semi-continuous and convex in the first variable and  for any natural number $l$, any points $(y_i,x_i)_{i=1}^l\in Y\times X$, $x\in X$ and any $\lambda_1,\dotsc,\lambda_l\geq 0$ such that $\sum_{i=1}^l \lambda_i=1$ we have
\begin{equation}\label{eqn:K-function}
\sum_{i,j=1}^l\lambda_i\lambda_j \Phi(y_i-y_j, x_i,x_j)\geq 2\sum_{i=1}^l \lambda_i \Phi(y_i-\sum_{j=1}^l\lambda_j y_j, x_i,x).
\end{equation}
If $Y$ is finite-dimensional, then it suffices to consider $l\leq \mathrm{dim}Y+1$.
\end{definition}

\begin{definition}
For a subset $S\subset Y$ of a real vector space $Y$ we denote by $\mathrm{Conv}S$ the smallest convex set in $Y$ containing $S$. It shall be called the convex hull of $S$.
\end{definition}

\begin{theorem}\label{thm:Minty}
Let $Y$ be a real vector space and $X$ be a set. Let  $\Phi\colon Y\times X\times X\to\mathbb{R}$ be a Kirszbraun function. Let 
\begin{equation*}
(y_i,x_i)_{i=1}^l\subseteq Y\times X
\end{equation*}
be a sequence such that
\begin{equation*}
\Phi(y_i-y_j,x_i,x_j)\leq 0
\end{equation*}
for all $i,j=1,\dotsc,l$.
Let $x\in X$. Then there exists a vector $y\in Y$ such that
\begin{equation*}
\Phi(y_i-y,x_i,x)\leq 0
\end{equation*}
for all $i=1,\dotsc,l$.
Furthermore, $y$ may be chosen to lie in $\mathrm{Conv}\{y_1,\dotsc,y_l\}$.
\end{theorem}

Let us mention that the proof of the above theorem relies on von Neumann's minimax theorem.

\section{Existence of extensions of $1$-Lipschitz maps}\label{s:existence}


\begin{proof}[Proof that \ref{i:v} implies \ref{i:K} in Theorem \ref{thm:main}]
Let us assume that \ref{i:v} holds true. We shall show that \ref{i:K} holds as well.
We begin with showing that for any natural number $k$, any points $x_1,\dotsc,x_k\in A$ and any $x\in X$ the intersection of closed, convex sets\footnote{$B(x,r)$ denotes a ball of radius $r>0$ centred at $x$.}
\begin{equation*}
\bigcap_{i=1}^k B(u(x_i),\norm{x_i-x})\cap ( v(x)+K)
\end{equation*}
is non-empty. 
By the Helly theorem, see \cite{Helly1923}, it is enough to prove that for any finite $m\leq\mathrm{dim}Y$ and any $k\leq m$ the above intersection is non-empty and that for $k\leq m+1$ the intersection 
\begin{equation*}
\bigcap_{i=1}^k B(u(x_i),\norm{x_i-x})
\end{equation*}
is non-empty. The latter follows by Kirszbraun's theorem, Theorem \ref{thm:Kirszbraun}. To prove the former, let
\begin{equation*}
    V=\mathrm{Aff}\{u(x_i)-v(x_i)\mid i=1,\dotsc,k\},
\end{equation*}
and consider a function $\Phi\colon V\times X\times X\to\mathbb{R}$ 
defined by the formula
\begin{equation*}
\Phi(y,x,x')=\norm{y}^2-2\langle y, v(x)-v(x')\rangle -\norm{x-x'}^2+\norm{v(x)-v(x')}^2, y\in V, x,x'\in X.
\end{equation*}
Let us check that it is a Kirszbraun function, see Definition \ref{defin:Ka}. The condition of convexity and finitely lower semi-continuity is clearly satisfied. We need only to check whether the condition (\ref{eqn:K-function}) holds, with $l\leq m$, as $\mathrm{dim}V\leq k-1\leq m-1$.  It is readily seen that the first two summands in the definition of $\Phi$ both satisfy the condition (\ref{eqn:K-function}) with equalities. Indeed, let $l\leq m$  and $(y_i,x_i)_{i=1}^l\in V\times X$, $x\in X$ and let $\lambda_1,\dotsc,\lambda_l\geq 0$ be such that $\sum_{i=1}^l\lambda_i=1$. Then a short calculation readily implies that
\begin{equation*}
\sum_{i.j=1}^l\lambda_i\lambda_j \norm{y_i-y_j}^2-2\sum_{i=1}^l\lambda_i \norm{y_i-\sum_{j=1}^l\lambda_jy_j}^2=0
\end{equation*}
and that
\begin{equation*}
\sum_{i.j=1}^l\lambda_i\lambda_j \langle y_i-y_j,x_i-x_j\rangle -2\sum_{i=1}^l\lambda_i \langle y_i-\sum_{j=1}^l\lambda_jy_j,x_i-x\rangle=0.
\end{equation*}
Thus, to satisfy (\ref{eqn:K-function}), we need to show that 
\begin{equation*}
\sum_{i,j=1}^l \lambda_i\lambda_j \big(\norm{v(x_i)-v(x_j)}^2-\norm{x_i-x_j}^2\big)\geq 2\sum_{i=1}^l \lambda_i \big(\norm{v(x_i)-v(x)}^2-\norm{x_i-x}^2\big).
\end{equation*}
Rearranging we get an equivalent inequality
\begin{equation*}
\Big\lVert \sum_{i=1}^l \lambda_ix_i-x\Big\rVert ^2\geq \Big\lVert \sum_{i=1}^l \lambda_iv(x_i)-v(x)\Big\rVert^2,
\end{equation*}
which holds true by the assumption on $v$.

By Theorem \ref{thm:Minty} we infer that for any points $x_1,\dotsc,x_k\in A$, $k\leq m$, and any $x\in X$ the intersection of closed, convex sets
\begin{equation}\label{eqn:inter}
\bigcap_{i=1}^k B(u(x_i),\norm{x_i-x})\cap ( v(x)+K)
\end{equation}
is non-empty. Indeed, let us put $y_i=u(x_i)-v(x_i)$ for $i=1,2,\dotsc,k$. Then
\begin{equation*}
\Phi(y_i-y_j,x_i,x_j)\leq 0\text{ for }i,j=1,2,\dotsc,k
\end{equation*}
is equivalent to $u\colon \{x_1,\dotsc,x_k\}\to Y$ being $1$-Lipschitz.
Now, Theorem \ref{thm:Minty} ensures that there exists $y\in Y$ such that  
\begin{equation*}
\Phi(y-y_i,x_i,x)\leq 0\text{, for }i=1,\dotsc,k,
\end{equation*}
and such that $y\in K$, since $y\in\mathrm{Conv}\{y_1,\dotsc,y_k\}$.
The definition of $\Phi$ shows that this $y$ satisfies
\begin{equation*}
\norm{y+v(x)-(y_i+v(x_i))}\leq\norm{x-x_i}\text{ for }i=1,2,\dotsc,k.
\end{equation*}
Hence, setting $u(x)=y+v(x)$, we see that it belongs to the set defined by formula (\ref{eqn:inter}).

Thanks to the Banach--Alaoglu theorem the balls in $Y$ are compact in the weak* topology. Therefore for $x\in X$, the intersection of the form
\begin{equation*}
    \bigcap \Big\{B(u(x'),\norm{x-x'})\cap (v(x)+K)\mid x'\in A\Big\}
\end{equation*}
is non-empty if and only if any of its finite subintersections are non-empty. The latter we have already shown. 
Therefore we may always find $u(x)$, so  that $u$ is a $1$-Lipschitz map on $A\cup\{x\}$ and that 
\begin{equation*}
u(x)-v(x)\in K.
\end{equation*}
Let us order by inclusion all subsets of $X$ containing $A$ that admit a desired extension. By the Kuratowski--Zorn lemma, there exists a maximal subset with respect to this ordering. If it were not $X$, then by the above considerations we could find a strictly larger subset with an extension that satisfies the desired conditions. 
This is to say, \ref{i:K} is proven.
\end{proof}

Clearly, \ref{i:K} implies that \ref{i:ball} holds true. What remains to show is that \ref{i:ball} implies \ref{i:v} if $\mathrm{dim}Y\in\{1,2,3\}$ or if $X$ is convex.

\begin{remark}\label{rem:extend}
Let us note that a statement similar to that of Theorem \ref{thm:main}, concerned with extensions of maps $u\colon A\to Y$ satisfying 
\begin{equation}\label{eqn:sumext}
\Big\lVert u(x)-\sum_{i=1}^m t_iu(x_i)\Big\rVert\leq \Big\lVert x-\sum_{i=1}^mt_ix_i\Big\rVert
\end{equation}
for all $x,x_1,\dotsc,x_m\in A$ and all non-negative real numbers $t_1,\dotsc,t_m$ that sum up to one,
to maps $\tilde{u}\colon X\to Y$ such that (\ref{eqn:sumext}) holds true with  $x,x_1,\dotsc,x_m\in X$ and with non-negative real numbers $t_1,\dotsc,t_m$ that sum up to one, is false in general. Indeed, let us take $X=[0,1]^2$,  $A$ to be the set of vertices of $X$, $Y=\mathbb{R}^3$. Then any $\tilde{u}\colon X\to Y$ satisfying (\ref{eqn:sumext}) with $x,x_1,x_2,x_3\in X$ is necessarily affine, by Proposition \ref{pro:affine}. However, there exists $u\colon A\to Y$ that satisfies (\ref{eqn:sumext}), where $x,x_1,x_2,x_3\in A$, for which
\begin{equation*}
   \frac12( u(a_1)+u(a_4))\neq \frac12(u(a_2)+u(a_3)),
\end{equation*}
where $a_1=(0,0),a_2=(0,1),a_3=(1,0),a_4=(1,1)$.
Indeed, let us set 
\begin{equation*}
    u(a_1)=u(a_2)=u(a_3)=0\text{ and }u(a_4)=\frac1{\sqrt2}w.
\end{equation*}
for a unit vector $w\in Y$. Then (\ref{eqn:sumext}) holds true with $x,x_1,x_2,x_3\in A$, as
\begin{equation*}
    \inf\Big\{\Big\lVert x-\sum_{i=1}^3t_ix_i\Big\rVert\mid x,x_1,x_2,x_3\in A, t_1,t_2,t_3\geq 0, \sum_{i=1}^3t_i=1\Big\}=\frac1{\sqrt2}.
\end{equation*}
Therefore, there is no extension $\tilde{u}\colon X\to Y$ of $u$ that would satisfy the required inequalities.
\end{remark}


\section{Necessity of the condition}\label{s:necess}

We shall provide lemmata that we shall use to prove that if $\mathrm{dim}Y\in\{1,2,3\}$, then \ref{i:ball} of Theorem \ref{thm:main} implies \ref{i:v} of that theorem. As before, $Y,Z$ denote Hilbert spaces.

\subsection{Leaves of $1$-Lipschitz maps}\label{ss:leaves}

We shall recall \cite[Lemma 2.4]{Ciosmak20212}. Although it is formulated in \cite{Ciosmak20212} merely for finite-dimensional spaces, the proof works mutatis mutandis also for general Hilbert spaces. 

Let $A\subset Z$. We shall say that a map $v\colon A\to Y$ is an isometry provided that for all $x,y\in A$ we have $\norm{v(x)-v(y)}=\norm{x-y}$.

\begin{lemma}\label{lemma:unique}
Let $\mathcal{S}\subset Z$ be an arbitrary subset. Let $u\colon \mathcal{S}\to Y$ be an isometry. Then there exists a unique $1$-Lipschitz function $\tilde{u}\colon \mathrm{Conv}(\mathcal{S})\to Y$ such that $\tilde{u}|_{\mathcal{S}}=u$. Moreover $\tilde{u}$ is an affine isometry.
\end{lemma}


We cite the following definition, see \cite[Definition 2.2]{Ciosmak20212}, that is motivated by the above lemma. We refer the reader to \cite{Ciosmak20212} for a detailed study of this notion.

\begin{definition}
Let $u\colon Z\to Y$ be a $1$-Lipschitz map. A set $\mathcal{S}\subset Y$ is called a \emph{leaf} of $u$ if $u|_{\mathcal{S}}$ is an isometry and for any $y\notin \mathcal{S}$ there exists $x\in\mathcal{S}$ such that $\norm{u(y)-u(x)}<\norm{y-x}$.
\end{definition}

In other words, $\mathcal{S}$ is a leaf if it is a maximal set, with respect to the order induced by inclusion, such that $u|_{\mathcal{S}}$ is an isometry.

\begin{lemma}\label{lem:dec}
 Let $Z,Y$ be Hilbert spaces. 
Let $x\in Z$. Suppose that $A\subset Z$ and that $u\colon A\cup\{x\}\to Y$ is $1$-Lipschitz and is an isometry on $A$, i.e., for all $x_1,x_2\in A$
\begin{equation*}
    \norm{u(x_1)-u(x_2)}=\norm{x_1-x_2}.
\end{equation*}
Then for all $x_1,\dotsc,x_m\in A$ and all non-negative $t_1,\dotsc,t_m$ that sum up to one we have
\begin{equation*}
    \Big\lVert u(x)-\sum_{i=1}^mt_iu(x_i)\Big\rVert\leq \Big\lVert x-\sum_{i=1}^mt_ix_i\Big\rVert.
\end{equation*}
\end{lemma}
\begin{proof}
   Let us take a $1$-Lipschitz extension $\tilde{u}$ of $u$ to $Z$, which exists thanks to the Kirszbraun theorem. By Lemma \ref{lemma:unique}, $\tilde{u}$ is an affine isometry on $\mathrm{Conv}A$. Therefore, for $x_1,\dotsc,x_m\in A$ and non-negative $t_1,\dotsc,t_m$ that sum up to one, we have
   \begin{equation*}
       \tilde{u}\Big(\sum_{i=1}^mt_ix_i\Big)=\sum_{i=1}^mt_iu(x_i).
   \end{equation*}
   Since $\tilde{u}$ is $1$-Lipschitz we have
   \begin{equation*}
   \Big\lVert\ u(x)-\tilde{u}\Big(\sum_{i=1}^mt_i x_i\Big)\Big\rVert\leq \Big\lVert x_{m+1}-\sum_{i=1}^mt_i x_i\Big\rVert.
   \end{equation*}
\end{proof}

\subsection{Consequence of existence of isometric embedding}\label{ss:conseq}

The strategy of the proof that \ref{i:ball} of Theorem \ref{thm:main} implies \ref{i:v} of that theorem that we adopt is the following. Let $v\colon X\to Y$. Let $x_1,\dotsc,x_m\in X$.
We shall construct an isometric embedding $u\colon\{x_1,\dotsc,x_m\}\to Y$
in such a way that:
\begin{enumerate}
    \item for any $i=1,2,\dotsc,m$, $u(x_i)$ is roughly at a distance $\delta$ to the corresponding vertex $v(x_i)$,
    \item the simplex $\mathrm{Conv}\{u(x_1),\dotsc,u(x_m)\}$ is shifted in direction of 
    \begin{equation*}
    v(x_{m+1})-\sum_{i=1}^mt_i v(x_i)
    \end{equation*}
     with respect to $\mathrm{Conv}\{v(x_1),\dotsc,v(x_m)\}$.
\end{enumerate} 
Then we can infer that \ref{i:v} of Theorem \ref{thm:main} is valid for $v$, provided that we pick $\delta$ sufficiently large. This conclusion is the matter of the following lemma.

\begin{lemma}\label{lemma:shift}
Let $Z,Y$ be two real Hilbert spaces. Suppose that $X\subset Z$ and that $v\colon X\to Y$ satisfies \ref{i:ball} of Theorem \ref{thm:main}.
Let  $m\leq\mathrm{dim}Y$ be a finite number. 
Let $A=\{x_1,\dotsc,x_m\}\subset X$. Let $x_{m+1}\in X$,  $t_1,\dotsc,t_m$ be non-negative numbers that sum up to one.   
Suppose that 
\begin{equation*}
    \Big\lVert v(x_{m+1})-\sum_{i=1}^{m}t_iv(x_i)\Big\rVert>\Big\lVert x_{m+1}-\sum_{i=1}^{m}t_ix_i\Big\rVert.
\end{equation*}
Let $C>0$. Suppose that $\delta$ is greater than
\begin{equation}\label{eqn:delta}
8(\mathrm{diam}A)^2+3C+ \Bigg(\frac{8\,(\mathrm{diam}A)^2+2C}{\Big\lVert v(x_{m+1})-\sum_{i=1}^{m}t_iv(x_i)\Big\rVert-\Big\lVert x_{m+1}-\sum_{i=1}^{m}t_ix_i\Big\rVert}\Bigg)^2.
\end{equation}
Then there exists no isometry $u\colon A\to Y$ such that 
\begin{equation*}
\abs{\norm{u(x)-v(x)}^2-\delta}\leq C\text{ for all }x\in A,
\end{equation*}
and such that for some non-negative numbers $s_1,\dotsc,s_m$ that sum up to one, some non-negative $\lambda$ we have
\begin{equation}\label{eqn:lambdaw}
\sum_{i=1}^m s_i(u(x_i)-v(x_i))=-\lambda\Big(v(x_{m+1})-\sum_{i=1}^mt_iv(x_i)\Big).
\end{equation}
\end{lemma}

Before we prove the lemma, let us first prove the following.

\begin{lemma}\label{lem:conv}
Let $Y,Z$ be Hilbert spaces, $A\subset Z$.    Suppose that $u,v\colon A\to Y$ are $1$-Lipschitz. Let us assume that for some $\delta>0$ and $C>0$ there is 
\begin{equation*}
    \abs{\norm{u(x)-v(x)}^2-\delta}\leq C\text{ for all }x\in A.
\end{equation*}
Then for any $x_1,\dotsc,x_k\in A$ and any non-negative $t_1,\dotsc,t_k$ that add up to one we have
\begin{equation*}
\bigg\lvert\Big\lVert \sum_{i=1}^k t_i(u(x_i)-v(x_i))\Big\rVert^2-\delta\bigg\lvert\leq 8\big(\mathrm{diam}A\big)^2+3C.
\end{equation*}
\end{lemma}
\begin{proof}
   Since 
\begin{equation*}
    \abs{\norm{u(x_i)-v(x_i)}^2-\delta}\leq C\text{ for }i=1,\dotsc,k,
\end{equation*}
we see that for $i,j=1,\dotsc,k$
\begin{equation*}
\Big\lvert2\Big\langle u(x_i)-v(x_i)-\big(u(x_j)-v(x_j)\big), u(x_j)-v(x_j)\Big\rangle+\norm{u(x_i)-v(x_i)-(u(x_j)-v(x_j))}^2\Big\rvert\leq 2C.
\end{equation*}
Indeed, after expanding the square, the above is equivalent to
\begin{equation*}
\big\lvert \norm{u(x_i)-v(x_i)}^2-\norm{u(x_j)-v(x_j)}^2\big\rvert\leq 2C,
\end{equation*}
and thus follows by the assumption and the triangle inequality.
We see now that for each $j=1,\dotsc,k$
\begin{align*}
&\Big\lvert \Big\lVert \sum_{i=1}^k t_i(u(x_i)-v(x_i))\Big\rVert^2-\delta\Big\rvert=\\
&\Big\lvert \Big\lVert \sum_{i=1}^k t_i\big(u(x_i)-v(x_i)-(u(x_j)-v(x_j))\big)+u(x_j)-v(x_j)\Big\rVert^2-\delta\Big\rvert=\\
&\Big\lvert\Big\lVert \sum_{i=1}^k t_i(u(x_i)-v(x_i)-(u(x_j)-v(x_j)))\Big\rVert^2+\big(\norm{u(x_j)-v(x_j)}^2-\delta\big)+\\
&2\sum_{i=1}^kt_i \langle u(x_i)-v(x_i)-(u(x_j)-v(x_j)), u(x_j)-v(x_j)\rangle\Big\rvert\leq\\
& \Big\lvert\Big\lVert \sum_{i=1}^k t_i(u(x_i)-v(x_i)-(u(x_j)-v(x_j)))\Big\rVert^2+\sum_{i=1}^kt_i \norm{ u(x_i)-v(x_i)-(u(x_j)-v(x_j))}^2\Big\rvert+3C.
\end{align*}
Thanks to the $1$-Lipschitzness of $u$ and of $v$ we infer that
\begin{equation*}
\bigg\lvert\Big\lVert \sum_{i=1}^k t_i(u(x_i)-v(x_i))\Big\rVert^2-\delta\bigg\lvert\leq 8\big(\mathrm{diam}A\big)^2+3C.
\end{equation*}
\end{proof}

\begin{proof}[Proof of Lemma \ref{lemma:shift}]
Let us suppose on the contrary, that an isometry $u\colon A\to Y$ with the outlined properties exists.
Since 
\begin{equation*}
    \abs{\norm{u(x_i)-v(x_i)}^2-\delta}\leq C\text{ for }i=1,\dotsc,m,
\end{equation*}
Lemma \ref{lem:conv} shows that
\begin{equation}\label{eqn:bound}
\bigg\lvert\Big\lVert \sum_{i=1}^m t_i(u(x_i)-v(x_i))\Big\rVert^2-\delta\bigg\lvert\leq 8\big(\mathrm{diam}A\big)^2+3C.
\end{equation}
Clearly, $u$ is $1$-Lipschitz on $A$. By the assumption on $v$ and by Theorem \ref{thm:main}, there exists a $1$-Lipschitz extension $\tilde{u}$ of $u$ to $A\cup \{x_{m+1}\}$ such that \begin{equation}\label{eqn:condi}
\norm{\tilde{u}(x_{m+1})-v(x_{m+1})}\leq\big(\delta+C\big)^{\frac12}.
\end{equation}
By Lemma \ref{lem:dec} 
\begin{equation}\label{eqn:uum}
\Big\lVert\tilde{u}(x_{m+1})-\sum_{i=1}^mt_i u(x_i)\Big\rVert\leq \Big\lVert x_{m+1}-\sum_{i=1}^mt_i x_i\Big\rVert.
\end{equation}
Let us set 
\begin{align*}
& \Delta_u=\tilde{u}(x_{m+1})-\sum_{i=1}^mt_i u(x_i),\\
& \Delta_v=v(x_{m+1})-\sum_{i=1}^m t_iv(x_i)\\
& \Delta_x=x_{m+1}-\sum_{i=1}^m t_ix_i.
\end{align*}
Condition (\ref{eqn:condi}) implies that
\begin{equation*}
\Big\lVert \Delta_u-\Delta_v+\sum_{i=1}^mt_i(u(x_i)-v(x_i))\Big\rVert\leq \big(\delta+C\big)^{\frac12}.
\end{equation*}
By the Pythagorean theorem we see that
\begin{equation*}
2\Big\langle \Delta_u-\Delta_v, \sum_{i=1}^m t_i(u(x_i)-v(x_i))\Big\rangle+\Big\lVert \sum_{i=1}^m t_i(u(x_i)-v(x_i))\Big\rVert^2+\norm{\Delta_u-\Delta_v}^2\leq \delta+C.
\end{equation*}
Therefore, by (\ref{eqn:bound}), we get that
\begin{equation*}
2\Big\langle \Delta_u-\Delta_v, \sum_{i=1}^m t_i(u(x_i)-v(x_i))\Big\rangle\leq C+\Big\lvert \Big\lVert \sum_{i=1}^m t_i(u(x_i)-v(x_i))\Big\rVert^2-\delta\Big\rvert\leq 8(\mathrm{diam}A)^2+4C.
\end{equation*}
By (\ref{eqn:lambdaw}), $2\Big\langle \Delta_u-\Delta_v, \sum_{i=1}^m t_i(u(x_i)-v(x_i))\Big\rangle$ is equal to
\begin{equation*}
-2\lambda\Big\langle \Delta_u-\Delta_v,\Delta_v\Big\rangle+2\Big\langle\Delta_u-\Delta_v ,\sum_{i,j=1}^m t_is_j(u(x_i)-v(x_i)-u(x_j)+v(x_j))\Big\rangle.
\end{equation*}
Thus, thanks to the $1$-Lipschitzness of $u$ and of $v$, we get
\begin{equation*}
-2\lambda\Big\langle \Delta_v,\Delta_u-\Delta_v\Big\rangle\leq 16(\mathrm{diam}A)^2+4C.
\end{equation*}
Since by  (\ref{eqn:uum}) $\norm{\Delta_u}\leq\norm{\Delta_x}$, we have
\begin{equation*}
2\lambda\norm{\Delta_v}^2\leq 16(\mathrm{diam}A)^2+4C+2\lambda \norm{\Delta_v}\norm{\Delta_x}.
\end{equation*}
By  the assumption, $\norm{\Delta_v}>\norm{\Delta_x}$. Then
\begin{equation*}
\lambda\norm{\Delta_v}\leq \frac{8(\mathrm{diam}A)^2+2C}{\norm{\Delta_v
}-\norm{\Delta_x}}.
\end{equation*}
 Similarly as in (\ref{eqn:bound}) we show that, using (\ref{eqn:lambdaw}),
\begin{equation*}
\abs{\lambda^2\norm{\Delta_v}^2-\delta}\leq 8(\mathrm{diam}A)^2+3C.
\end{equation*}
It follows that
\begin{equation*}
\delta\leq \lambda^2\norm{\Delta_v}^2+8(\mathrm{diam}A)^2+3C\leq 8(\mathrm{diam}A)^2+3C+ \bigg(\frac{8\,(\mathrm{diam}A)^2+2C}{\norm{\Delta_v}-\norm{\Delta_x}}\bigg)^2.
\end{equation*}
This stands in contradiction with the assumption that $\delta$ is greater than (\ref{eqn:delta}).
\end{proof}

\begin{remark}\label{rem:relax}
    Let us observe that the assumption that $u\colon A\to Y$ is an isometry was used merely to infer, by a use of Lemma \ref{lem:dec}, that for any $1$-Lipschitz extension $\tilde{u}\colon A\cup\{x_{m+1}\}\to Y$ of $u$, we have
    \begin{equation*}
        \Big\lVert \tilde{u}(x_{m+1})-\sum_{i=1}^mt_iu(x_i)\Big\rVert\leq  \Big\lVert x_{m+1}-\sum_{i=1}^mt_ix_i\Big\rVert,
    \end{equation*}
    for all $x_1,\dotsc,x_m\in A$ and all non-negative $t_1,\dotsc,t_m$ that sum up to one.
    Thus, the same conclusion can be drawn from the existence of $u$ that satisfies the above inequality for any of its $1$-Lipschitz extension to $A\cup \{x\}$.
\end{remark}

\begin{remark}
When we assume that $v$ satisfies \ref{i:K} of Theorem \ref{thm:main}, then the proof of an analogue of Lemma \ref{lem:conv} is considerably simpler. The essence of the idea is however the same. 
\end{remark}

\subsection{Existence of embedding and completion of the proof for $\mathrm{dim}Y\leq 3$}\label{ss:construction}

We shall now provide a construction of a map $u$, as in Lemma \ref{lemma:shift}, for $m\leq \min \{\mathrm{dim}Y,3\}$. 

\begin{lemma}\label{lem:construction}
Let $Z,Y$ be two real Hilbert spaces. Suppose that $X\subset Z$ and that $v\colon X\to Y$ satisfies \ref{i:ball} of Theorem \ref{thm:main}.
Let  $m\leq \min \{\mathrm{dim}Y,3\}$, $A=\{x_1,\dotsc,x_m\}\subset X$. Then for any unit vector $w\in Y$, $\delta>0$
 there exists an isometry $u\colon A\to Y$ such that 
\begin{equation}\label{eqn:C}
\abs{\norm{u(x)-v(x)}^2-\delta}\leq 4(\mathrm{diam}A)^2\text{ for all }x\in A,
\end{equation}
and such that 
\begin{equation*}
u(x_1)-v(x_1)=\sqrt{\delta}w.
\end{equation*}
\end{lemma}
\begin{proof}
We shall construct an isometry $u\colon A\to Y$ such that $u(x_1)-v(x_1)=\sqrt{\delta}w$ and such that
\begin{equation}\label{eqn:dew}
    \langle u(x)-v(x),w\rangle=\sqrt{\delta}\text{ for all }x\in A.
\end{equation}
It will follow that (\ref{eqn:C}) is also satisfied.
Indeed, (\ref{eqn:dew}) implies trivially that
\begin{equation*}
    \norm{u(x)-v(x)}^2\geq \delta\text{ for all }x\in A.
\end{equation*}
Moreover
\begin{align*}
    & \norm{u(x)-v(x)}^2= \norm{u(x)-u(x_1)-(v(x)-v(x_1))+\sqrt{\delta}w}^2=\\
    &-\delta+2\sqrt{\delta}\langle u(x)-v(x),w\rangle + \norm{u(x)-u(x_1)-(v(x)-v(x_1))}^2\leq \delta+4(\mathrm{diam}A)^2.
\end{align*}
Thus indeed (\ref{eqn:C}) will follow. 

If $m=1$, then we define 
\begin{equation*}
 u(x_1)=v(x_1)+\sqrt{\delta}w.
\end{equation*}
Then  $u$ clearly satisfies our requirements.

Let $m=2$. Clearly, we may assume that $x_1\neq x_2$. We set
\begin{equation*}
 u(x_1)=v(x_1)+\sqrt{\delta}w.
\end{equation*}
We shall now find $u(x_2)\in Y$ such that 
\begin{equation}\label{eqn:req}
    \norm{u(x_1)-u(x_2)}=\norm{x_1-x_2}\text{ and }\langle u(x_2)-v(x_2),w\rangle^2=\delta.
\end{equation}
We may look for $u(x_2)$ in the affine space spanned by $\{v(x_1),v(x_2),u(x_1)\}$, whose dimension is at most two.
Consider the line spanned by $v(x_1)$ and $u(x_1)$. It suffices to show that  
\begin{equation*}
    y_1=u(x_1)+\norm{x_1-x_2}w,y_2=u(x_1)-\norm{x_1-x_2}w
\end{equation*}
satisfy
\begin{equation}\label{eqn:ul}
    \langle y_1-v(x_2),w\rangle\geq\sqrt{\delta},    \langle y_2-v(x_2),w\rangle\leq\sqrt{\delta}.
\end{equation}
Note that by Lemma \ref{lemma:shift}, and the construction for $m=1$, it follows that $v$ is $1$-Lipschitz. Therefore, we get
\begin{equation*}
    \langle y_1-v(x_2),w\rangle\geq \norm{x_1-x_2}+\langle u(x_1)-v(x_1),w\rangle -\langle v(x_2)-v(x_1),w\rangle \geq\sqrt{\delta}.
\end{equation*}
Similarly we show the other inequality. Now, the possible positions of $u(x_2)$, for which $u\colon \{x_1,x_2\}\to Y$ is an isometry, form a circle in the space spanned by $\{u(x_1),v(x_1),v(x_2)\}$. Thanks to connectedness of the circle and thanks to (\ref{eqn:ul}), we infer the existence of $u(x_2)$ that satisfies (\ref{eqn:req}).


Let now $m=3$. We may assume that $\mathrm{dim}Y=3$. 

For any distinct numbers $i,j\in\{1,2,3\}$ we may find $u(x_i),u(x_j)$ in the way described above. We infer by Lemma \ref{lemma:shift}
 that for any $j\in \{1,2,3\}$ and any non-negative real numbers $(r_i)_{i\in \{1,2,3\}\setminus \{j\}}$  that sum up to one we have
 \begin{equation*}
     \Big\lVert v(x_j)-\sum_{i\in \{1,2,3\}\setminus\{j\}}r_iv(x_i)\Big\rVert\leq \Big\lVert x_j-\sum_{i\in \{1,2,3\}\setminus\{j\}}r_ix_i\Big\rVert.
 \end{equation*}
 This is equivalent to
  \begin{equation}\label{eqn:vv}
     \Big\lVert \sum_{i=1}^3s_iv(x_i)\Big\rVert\leq \Big\lVert \sum_{i=1}^3s_ix_i\Big\rVert\text{ for all }s_1,s_2,s_3\in\mathbb{R}\text{ such that }\sum_{i=1}^3s_i=0.
 \end{equation}
We set 
\begin{equation}\label{eqn:1}
    u(x_1)=v(x_1)+\sqrt{\delta}w
\end{equation}
and construct $u(x_2)$, lying in the two-dimensional 
affine space $W$ containing the span of $u(x_1),v(x_1),v(x_2)$, such that 
\begin{equation}\label{eqn:2}
    \norm{u(x_1)-u(x_2)}=\norm{x_1-x_2}\text{ and }\langle u(x_2)-v(x_2),w\rangle^2=\delta.
\end{equation}

Next, we shall construct $u(x_3)$, so that $u\colon\{x_1,x_2,x_3\}\to Y$ is an isometry and such that (\ref{eqn:dew}) holds true.
Observe that the possible positions for $u(x_3)$, for which $u\colon \{x_1,x_2,x_3\}\to Y$ is an isometry, lie on a circle contained in $Y$. 
The intersection of the considered circle with $W$ consists of two points $y_+,y_-$. 
There are two possible positions of $v(x_3)$ in $W$, given the positions of $v(x_1),v(x_2)$, which we shall denote by $z_+$ and $z_-$. 
The situation is illustrated on Figure \ref{fig:triangles}.

\begin{figure}
\begin{tikzpicture}
\draw[thick] (4,0) -- (6,0) -- (5.7,2) -- cycle;
\draw[thick] (4,0) --  (5.7,-2)--(6,0);
\coordinate (A) at (4,0);
\coordinate (B) at (6,0);
\coordinate (C) at (5.7,-2);
\coordinate (D) at (5.7,2);
\draw[fill] (A) circle (1.0pt);
\draw[fill] (B) circle (1.0pt);
\draw[fill] (C) circle (1.0pt);
\draw[fill] (D) circle (1.0pt);
\node[xshift=-14pt] at (A) {$u(x_1)$};
\node[xshift=14pt] at (B) {$u(x_2)$};
\node[yshift=7pt] at (D) {$y_+$};
\node[yshift=-7pt] at (C) {$y_-$};

\draw[thick] (9,-4) -- (11,-4) -- (10.35,-2.8) -- cycle;
\draw[thick] (9,-4)  -- (10.35,-5.2)--(11,-4) ;
\coordinate (a) at (9,-4);
\coordinate (b) at (11,-4);
\coordinate (c) at (10.35,-2.8);
\coordinate (d) at (10.35,-5.2);
\draw[fill] (a) circle (1.0pt);
\draw[fill] (b) circle (1.0pt);
\draw[fill] (c) circle (1.0pt);
\draw[fill] (d) circle (1.0pt);
\node[xshift=-14pt] at (a) {$v(x_1)$};
\node[xshift=14pt] at (b) {$v(x_2)$};
\node[yshift=7pt] at (c) {$z_+$};
\node[yshift=-7pt] at (d) {$z_-$};

\draw[->,thick] (12,-1.5)--(9.5,0.5);
\node at (10.75,-1) {$w$};
\end{tikzpicture}
\caption{The illustration of exemplary positions of $u(x_1),u(x_2)$, $v(x_1),v(x_2)$ and of $y_-,y_+,z_-,z_+$ on the $W$ plane.} \label{fig:triangles}
\end{figure}
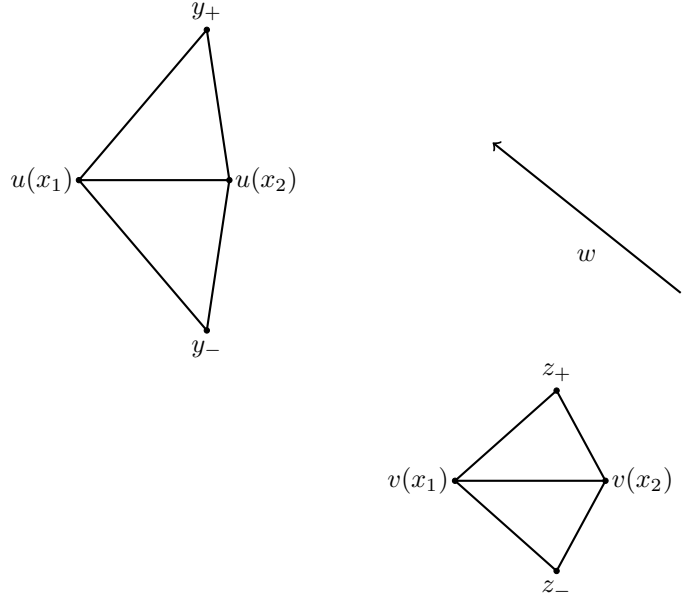

Note that $W$ is equal to  the affine hull of $u(x_1),u(x_2),y_+$, and is equal to the affine hull of $u(x_1),u(x_2),y_-$.
Therefore
\begin{equation}\label{eqn:ww}
    w=t_1u(x_1)+t_2u(x_2)+t_3y_+
\end{equation}
 for some $t_1,t_2,t_3\in\mathbb{R}$ that sum up to zero with $t_3\geq 0$ and 
 \begin{equation}\label{eqn:ww2}
    w=s_1u(x_1)+s_2u(x_2)+s_3y_{-}
\end{equation}
 for some $s_1,s_2,s_3\in\mathbb{R}$ that sum up to zero with $s_3\leq 0$. 
 The coefficients $t_3,s_3$ are non-zero, unless $w$ is proportional to $u(x_1)-u(x_2)$. Let us assume that it is not.
We shall show that 
\begin{equation}\label{eqn:circle}
    \langle P_U(y_+-v(x_3)),w\rangle\geq\sqrt{\delta},  \langle P_U(y_--v(x_3)),w\rangle\leq\sqrt{\delta},
\end{equation}
where $P_U$ is the orthogonal projection onto the tangent space $U$ to $W$. Hence, we may assume that $v(x_3)$ belongs to $W$, and thus it suffices to prove inequalities  (\ref{eqn:circle}) with $v(x_3)$ replaced by $z_+$ and by $z_-$ respectively.

By (\ref{eqn:vv}), the fact that $u$ extends to an unique isometry on the affine hull of $\{x_1,x_2,x_3\}$, see Lemma \ref{lemma:unique}, and the Cauchy--Schwarz inequality, we get
\begin{equation*}
    \big\lVert t_1u(x_1)+t_2u(x_2)+t_3y_+\big\rVert^2\geq \big\langle t_1v(x_1)+t_2v(x_2)+t_3z_+,t_1u(x_1)+t_2u(x_2)+t_3y_+\big\rangle,
\end{equation*}
which implies, together with (\ref{eqn:ww}), that
\begin{equation*}
    \big\langle  \big(t_1u(x_1)+t_2u(x_2)+t_3y_+\big)-\big(t_1v(x_1)+t_2v(x_2)+t_3z_+\big),w\big\rangle\geq  0,
\end{equation*}
and consequently, thanks to (\ref{eqn:1}) and to (\ref{eqn:2}), that
\begin{equation}\label{eqn:t}
   t_3 \langle y_+-z_+,w\rangle\geq  -\sum_{i=1}^2t_i\langle  u(x_i)-v(x_i),w\rangle= t_3\sqrt{\delta}.
\end{equation}
Similarly, 
\begin{equation}\label{eqn:s}
    s_3\langle y_--z_-,w\rangle\geq - \sum_{i=1}^2s_i\langle  u(x_i)-v(x_i),w\rangle=s_3\sqrt{\delta}. 
\end{equation}
Since $t_3>0$ and $s_3<0$, this shows that (\ref{eqn:circle}) holds true. Since the circle is connected, this implies the existence of $u(x_3)$ that satisfies our requirements.

 Let us now suppose $w$ and $u(x_1)-u(x_2)$ are proportional. Then $t_3$ and $s_3$ are zero. Then, in the above inequalities (\ref{eqn:t}), (\ref{eqn:s}) we get equalities. Therefore, we also have equality in the Cauchy--Schwarz inequality, which implies that there is $\lambda>0$ such that \begin{equation*}
  \lambda(   v(x_1)-v(x_2))=u(x_1)-u(x_2). 
 \end{equation*}
Using (\ref{eqn:2}), we see that $\norm{v(x_1)-v(x_2)}=\norm{x_1-x_2}$. 
Thus, in (\ref{eqn:circle}), we get equalities, which also shows that $u(x_3)$ satisfying our requirements exists.

Lemma \ref{lemma:shift} completes the proof.
\end{proof}

\begin{proof}[Proof that \ref{i:ball} implies \ref{i:v} in Theorem \ref{thm:main}, when $\mathrm{dim}Y\leq 3$]
    The proof is a direct consequence of Lemma \ref{lemma:shift} and Lemma \ref{lem:construction}.
\end{proof}
\begin{remark}
    Let us observe that the proof shows that \ref{i:ball} in Theorem \ref{thm:main} implies \ref{i:v} for all $m\leq\min\{\mathrm{dim}Y,3\}$.
\end{remark}

\begin{example}
    Let us note that one may find $v\colon X\to Y$, with $\mathrm{dim}Y=4$, such that \ref{i:v} of Theorem \ref{thm:main} is satisfied, yet there exist $A=\{x_1,x_2,x_3,x_4\}\subset X$, $\delta>0, C>0$ with no isometric embedding $u\colon A\to Y$ satisfying 
    \begin{equation}\label{eqn:w}
        \abs{\langle u(x)-v(x),w\rangle -\delta}\leq C\text{ for all }x\in A
    \end{equation}
    and
    $u(x_1)-v(x_1)=\delta w$ for a unit vector $w\in Y$. Indeed, let us take $X$ to consist of vertices $\{(0,0),(0,1),(1,0),(1,1)\}$ of $[0,1]^2$. Let 
    \begin{equation*}
        v((0,0))=v((0,1))=v((1,0))=0\text{ and let }v((1,1))=\frac1{\sqrt{2}}w.
    \end{equation*} 
    Then $v$ satisfies \ref{i:v} of Theorem \ref{thm:main}, cf. Remark \ref{rem:extend}.
    We shall specify $C>0$ and $\delta>0$ later.
    Suppose that we can find an isometry $u\colon A\to Y$ with the required properties. Then \begin{equation}\label{eqn:op}
        \abs{\langle u(x),w\rangle-\delta}\leq C\text{ for all }x\in \{(0,0),(0,1),(1,0)\}.
    \end{equation}
    Moreover, Lemma \ref{lemma:unique} and the Kirszbraun theorem show that 
    \begin{equation}\label{eqn:uo}
       u((1,1))= u((0,1))+u((1,0))-u((0,0)).
    \end{equation}
    Therefore
\begin{equation*} 
   \abs{ \langle u((1,1)),w\rangle-\delta}\leq 3C.
\end{equation*}
Note that 
\begin{equation*}
    u((1,1))-v((1,1))=u((1,1))-\frac1{\sqrt{2}}w = u((1,1))-u((0,0))+\Big(\delta-\frac1{\sqrt{2}}\Big)w.
\end{equation*}
It follows that
\begin{equation*}
     \abs{\langle u((1,1))-v((1,1)),w \rangle-\delta}=  \Big\lvert\langle u((1,1))-u((0,0)),w \rangle-\frac1{\sqrt{2}}\Big\rvert.
\end{equation*}
By the assumption, the above should be at most $C$. On the other hand, by (\ref{eqn:uo}) and by (\ref{eqn:op}), we see that it is bounded from below by $\frac{1}{\sqrt{2}}-4C$.
If we thus pick $C$ so that 
\begin{equation*}
    C<\frac{1}{5\sqrt{2}},
\end{equation*}
then we get a contradiction.

Similarly, we could show that there is no isometric embedding $u\colon A\to Y$ such that 
    \begin{equation*}
        \abs{\norm{ u(x)-v(x)}^2-\delta}\leq C\text{ for all }x\in A
    \end{equation*}
    in place of (\ref{eqn:w}).

The example shows that the idea, presented in this paper, of the proof that \ref{i:ball}, or \ref{i:K}, implies \ref{i:v} of Theorem \ref{thm:main}, does not work when $m>3$. 
\end{example}






\subsection{Case of convex $X$.}\label{ss:convex}

Let us now discuss the case when $X\subset Z$ is a convex set. To this aim we shall need the following observation.

\begin{proposition}\label{pro:affine}
    Let us suppose that $X$ is a convex set and $\mathrm{dim}Y\geq 2$. Let $v\colon X\to Y$. Then $v$ is $1$-Lipschitz and affine if and only if for all finite $m\leq \mathrm{dim}Y$, all $x,x_1,\dotsc,x_m\in X$ and all non-negative real numbers $t_1,\dotsc,t_m$ that sum up to one there is 
\begin{equation}\label{eqn:v}
\Big\lVert v(x)-\sum_{i=1}^mt_iv(x_i)\Big\rVert\leq \Big\lVert x-\sum_{i=1}^mt_ix_i\Big\rVert.
\end{equation}
\end{proposition}
\begin{proof}
  Let us   suppose that (\ref{eqn:v}) holds true and let $x_1,x_2\in X$, $t\in (0,1)$. Setting $x=tx_1+(1-t)x_2$, 
(\ref{eqn:v}) shows that $v(x)=tv(x_1)+(1-t)v(x_1)$. This is to say, $v$ is affine. Similarly, applying (\ref{eqn:v}) for $m=1$, we see that $v$ is $1$-Lipschitz.

Conversely, if $v$ is $1$-Lipschitz and affine and $m\leq\mathrm{dim}Y$, $x_1,\dotsc,x_m\in X$, $t_1,\dotsc,t_m$ are non-negative and sum up to one, then $\sum_{i=1}^m t_ix_i\in X$ and
\begin{equation*}
\sum_{i=1}^mt_iv(x_i)=v\Big(\sum_{i=1}^m t_ix_i\Big).
\end{equation*}
Thus, (\ref{eqn:v}) is satisfied thanks to the assumption that $v$ is $1$-Lipschitz.
\end{proof}

\begin{proof}[Proof that \ref{i:ball} implies \ref{i:v} in Theorem \ref{thm:main}, when $X$ is convex]
    Let us suppose that $X$ is convex and that \ref{i:ball} holds true. By Lemma \ref{lemma:shift} and Lemma \ref{lem:construction} we see that for $m\leq \min\{2,\mathrm{dim}Y\}$ and all $x,x_1,\dotsc,x_m\in X$ and all non-negative real numbers $t_1,\dotsc,t_m$ that sum up to one we have
\begin{equation*}
\Big\lVert v(x)-\sum_{i=1}^mt_iv(x_i)\Big\rVert\leq \Big\lVert x-\sum_{i=1}^mt_ix_i\Big\rVert.
\end{equation*}
Now, if $\mathrm{dim}Y\geq 3$, Proposition \ref{pro:affine} shows that the above inequality implies \ref{i:v} in general. 
\end{proof}

\begin{remark}
  As we have seen in Lemma \ref{lem:construction} the construction of a map $u$ that shows that \ref{i:ball} implies \ref{i:v} in Theorem \ref{thm:main} is rather simple for $m=2$. Our developments provide thus a more transparent proof of the main result of \cite[Theorem 3.9]{Ciosmak2021}, concerned with the case when $X$ is the entire Hilbert space.
\end{remark}

\section{Existence of extensions of monotone maps}\label{s:monotone}

The proof of Theorem \ref{thm:monotone} follows similar steps to that of Theorem \ref{thm:main}, thus we give a concise presentation.

\begin{proof}[Proof of Theorem \ref{thm:monotone}]
Let us assume that $v$ satisfies the assumption.
We begin with showing that for any points $x_1,\dotsc,x_k\in A$ and any $x\in X$ the intersection of closed, convex sets
\begin{equation}\label{eqn:inte}
\bigcap_{i=1}^k \Big\{y\in Y^*\mid \langle y-u(x_i),x-x_i\rangle\geq 0 \Big\}\cap ( v(x)+K)
\end{equation}
is non-empty. 
By the Helly theorem, see \cite{Helly1923}, it is enough to prove that for any finite $m\leq\mathrm{dim}Y$ and any $k\leq m$ the above intersection is non-empty and that for $k\leq m+1$ the intersection 
\begin{equation}\label{eqn:inters}
\bigcap_{i=1}^k \Big\{y\in Y^*\mid \langle y-u(x_i),x-x_i\rangle\geq  0\Big\}
\end{equation}
is non-empty. To prove the former,\footnote{The non-emptiness of (\ref{eqn:inters}) will follow by taking $v=0$ from the considerations for general $v$.} let
\begin{equation*}
    V=\mathrm{Aff}\{u(x_i)-v(x_i)\mid i=1,\dotsc,k\},
\end{equation*}
and consider a function $\Phi\colon V\times X\times X\to\mathbb{R}$ 
defined by the formula
\begin{equation*}
\Phi(y,x,x')=-\langle y+ v(x)-v(x'),x-x'\rangle, y\in V, x,x'\in X.
\end{equation*}
We need to check that it is a Kirszbraun function, see Definition \ref{defin:Ka}. The condition of convexity and finitely lower semi-continuity is clearly satisfied. We need only to check whether the condition (\ref{eqn:K-function}) holds, with $l\leq m$, as $\mathrm{dim}V\leq k-1\leq m-1$. Let $l\leq m$  and $(y_i,x_i)_{i=1}^l\in V\times X$, $x\in X$ and let $\lambda_1,\dotsc,\lambda_l\geq 0$ be such that $\sum_{i=1}^l\lambda_i=1$. Then a short calculation readily implies that (\ref{eqn:K-function}) for $\Phi$ is equivalent to
\begin{equation*}
\Big\langle v(x)-\sum_{i=1}^l \lambda_iv(x_i), x-\sum_{i=1}^l \lambda_ix_i\Big\rangle\geq 0,
\end{equation*}
which holds true by the assumption on $v$.

By Theorem \ref{thm:Minty} we infer that for any points $x_1,\dotsc,x_k\in A$, $k\leq m$, and any $x\in X$ the intersection (\ref{eqn:inte}) of closed, convex sets
is non-empty. Indeed, let
\begin{equation*}
y_i=u(x_i)-v(x_i)\text{ for }i=1,2,\dotsc,k.
\end{equation*}
Then the assumption that $u\colon \{x_1,\dotsc,x_k\}\to Y^*$ is monotone is equivalent to
\begin{equation*}
\Phi(y_i-y_j,x_i,x_j)\leq 0\text{ for }i,j=1,2,\dotsc,k.
\end{equation*}
Theorem \ref{thm:Minty} implies the existence of $y\in \mathrm{Conv}\{y_1,\dotsc,y_k\}\subset K$ such that
\begin{equation*}
\Phi(y-y_i,x_i,x)\leq 0\text{ for }i=1,2,\dotsc,k,
\end{equation*}
which is equivalent to
\begin{equation*}
\langle y+v(x)-(y_i+v(x_i)), x-x_i\rangle\geq 0\text{ for }i=1,2,\dotsc,k.
\end{equation*}
Hence, setting $u(x)=y+v(x)$, we see that it belongs to the set defined by formula (\ref{eqn:inte}).
Note that the set $v(x)+K$ is closed and bounded. By the Hahn--Banach theorem it is weakly closed and thanks to the Banach--Alaoglu theorem it is compact. Therefore for $x\in X$, the intersection of the form
\begin{equation*}
    \bigcap \Big\{\big\{y\in Y^*\mid \langle y-u(x'),x-x'\rangle\geq 0 \big\}\mid x'\in A\Big\}\cap ( v(x)+K)
\end{equation*}
is non-empty if and only if any of its finite subintersections are non-empty. The latter we have already shown. 
Therefore we may always find $u(x)$, so  that $u$ is a monotone map on $A\cup\{x\}$ and that 
\begin{equation*}
u(x)-v(x)\in K.
\end{equation*}
We finish the proof in the same way as we did in the proof of Theorem \ref{thm:main}.
\end{proof}

\section{Further generalisations}\label{s:further}

We mention that part of Theorem \ref{thm:main}, that \ref{i:v} implies \ref{i:K}, may be further generalised, in several possible directions. The proof is very similar to the one of Theorem \ref{thm:main}, and thus we do not present its details.

One of the possibilities concerns extensions of maps $u\colon A\to Y$ that satisfy
\begin{equation*}
    \Big\lVert \sum_{i=1}^k t_iu(x_i)\Big\rVert\leq   \Big\lVert \sum_{i=1}^k t_ix_i\Big\rVert,
\end{equation*}
for all $t_1,\dotsc,t_k\in\mathbb{R}$ with $\sum_{i=1}^kt_i=0$ and all $x_1,\dotsc,x_k\in A$. Here $k$ is a fixed natural number.\footnote{If we allow for varying, arbitrary large, number of points, then any $u$ that would satisfy thus defined conditions, would necessarily be a restriction of an affine function on the underlying space.} In this case an analogue of condition \ref{i:v} of Theorem \ref{thm:main} for a map $v\colon X\to Y$ would be 
\begin{equation*}
    \Big\lVert \sum_{i=1}^{(m+1)(k-1)} t_iv(x_i)\Big\rVert\leq   \Big\lVert \sum_{i=1}^{(m+1)(k-1)} t_ix_i\Big\rVert,
\end{equation*}
for all $t_1,\dotsc,t_{(m+1)(k-1)}\in\mathbb{R}$ with $\sum_{i=1}^{(m+1)(k-1)}t_i=0$, all $x_1,\dotsc,x_{(m+1)(k-1)}\in X$, with $m\leq \mathrm{dim}Y$. Let us note that when $(m+1)(k-1)\geq 2(\mathrm{dim}\,\mathrm{span}X+1)$, then any  $v$ that satisfies this condition, extends uniquely to an affine map on $\mathrm{span}X$. 

\bibliographystyle{siam}
\bibliography{biblio}

\end{document}